\documentclass[12pt,leqno]{article}

\usepackage{amsfonts}
\usepackage{amssymb}
\usepackage[all]{xy}

\newcommand{\cM}{{\cal M}}

\newcommand{\coker}{{\rm coker}}

\newcommand{\codiv}{{\rm codiv}}

\newcommand{\Hom}{{\rm Hom}}
\newcommand{\Ext}{{\rm Ext}}

\newcommand{\id}{{\rm id}}

\newcommand{\odd}{{\rm odd}}
\newcommand{\even}{{\rm even}}
\newcommand{\Cone}{{\rm Cone}}
\newcommand{\Gr}{{\rm Gr}}
\newcommand{\rel}{{\rm rel}}
\newcommand{\et}{{\rm et}}

\newcommand{\ZZ}{{\mathbb Z}}

\newcommand{\GG}{{\mathbb G}}
\newcommand{\RR}{{\mathbb R}}
\newcommand{\NN}{{\mathbb N}}

\newcommand{\QQ}{{\mathbb Q}}

\newcommand{\ra}{\rightarrow}
\def\rightepi{{\longrightarrow \kern-0.7em \rightarrow}}
\newcommand{\oplusm}{\mathop{\oplus}\limits}
\newcommand{\limm}{\mathop{\lim}\limits}

\newcommand{\notteilt}{{\,\not{\kern-0.075em|}\,}}
\def\antiddots{\mathinner{\mkern1mu\raise1pt\vbox{\kern7pt\hbox{.}}\mkern2mu
    \raise4pt\hbox{.}\mkern2mu\raise7pt\hbox{.}\mkern1mu}}

\begin{document}

\vspace*{15ex}

\begin{center}
{\LARGE\bf Refined and $l$-adic Euler Characteristics of \\
\vspace*{0.4cm} Nearly Perfect Complexes}\\
\bigskip
by\\
\bigskip
{\sc David Burns, Bernhard K\"ock} and {\sc Victor Snaith}
\end{center}

\bigskip

\begin{quote}
{\footnotesize {\bf Abstract}. We lift the Euler characteristic of
a nearly perfect complex to a relative algebraic $K$-group by
passing to its $l$-adic Euler characteristics.

}

\end{quote}

\section*{Introduction}

Let $X$ be a smooth projective, geometrically connected surface
over a finite field $k$ and assume that the Brauer group
$H^2_{{\rm et}}(X,\GG_m)$ is finite (as is widely believed to be
true). Let $G$ be a finite group acting on $X$ in such a way that
the canonical projection $\pi: X\ra X/G$ is \'etale. The
assumption that $\pi$ is \'etale ensures that $R\Gamma_{\rm
et}(X,\GG_m)$ can be represented by a bounded complex of
cohomologically trivial $\ZZ[G]$-modules $C^*$ (cf.\
Proposition~3.2 in \cite{CKPS}).

The following results may be found in \cite{Li1} and \cite{Li2}
(see also \S 3 in \cite{CKPS}):\\
(i) $H^i_{{\rm et}}(X,\GG_m) = 0$ for $i \ge 5$.\\
(ii) The groups $H^0_\et(X,\GG_m)$ and $H^4_\et(X,\GG_m)$ are
always finite and $H^2_\et(X,\GG_m)$ is finite by assumption.\\
(iii) The group $H^1_\et(X,\GG_m) = {\rm Pic}(X)$ is finitely
generated.\\
(iv) There is a canonical isomorphism $H^3_\et(X,\GG_m) \,\,
\tilde{\ra} \,\, \Hom(H^1_\et(X,\GG_m),\QQ/\ZZ)$.

We let $L_3$ denote the quotient of $H^1_\et(X,\GG_m)$ by its
torsion subgroup, and set $L_i := 0$ for $i\not= 3$. We let
$\tau_i$ denote the obvious isomorphism from $\Hom(L_i,\QQ/\ZZ)$
to the maximal divisible subgroup $H^i(C^*)_{{\rm div}}$ of
$H^i(C^*)$ (for any $i$). Then the triple $(C^*, (L_i)_{i \in
\ZZ}, (\tau_i)_{i \in \ZZ})$ constitutes a `nearly perfect
complex' in the sense of \cite{CKPS}. That is, $C^*$ is a bounded
complex of cohomologically trivial $\ZZ[G]$-modules and, in each
degree $i$, the quotient $H^i(C^*)_\codiv :=
H^i(C^*)/H^i(C^*)_{{\rm div}}$ is finitely generated, the module
$L_i$ is finitely generated and $\ZZ$-free, and the map
\[\tau_i:
\Hom_\ZZ(L_i,\QQ/\ZZ) \,\, \tilde{\ra} \,\, H^i(C^*)_{{\rm div}}\]
is a $\ZZ[G]$-equivariant isomorphism.

We recall that the main algebraic result of the paper \cite{CKPS}
is that to any nearly perfect complex $(C^*, (L_i)_i, (\tau_i)_i)$
 as above one can associate a canonical {\em Euler characteristic}
element $\chi(C^*)$ in the Grothendieck group $K_0(\ZZ[G])$ of all
finitely generated projective $\ZZ[G]$-modules which has image in
the Grothendieck group $G_0(\ZZ[G])$ of all finitely generated
$\ZZ[G]$-modules equal to
\[\chi_{{\rm coh}}(C^*):= \left([H^+(C^*)_\codiv] +
[\Hom_\ZZ(L_-,\ZZ)] \right) -  \left([H^-(C^*)_\codiv] +
[\Hom_\ZZ(L_+,\ZZ)] \right).\] (Here, we write $H^+(C^*)_\codiv$
for $\oplus_{i\, \even} H^i(C^*)_\codiv$, $L_-$ for $\oplus_{i\,
\odd} L_i$, and so on.)

We further recall that the intersection pairing on ${\rm CH}^1(X)
= {\rm Pic}(X) = H^1_\et(X,\GG_m) = H^1(C^*)$ induces a
$\ZZ[G]$-equivariant homomorphism
\[\mu: H^1(C^*)_\codiv = H^1(C^*) \ra \Hom(H^1(C^*), \ZZ) = \Hom(L_3,\ZZ)\]
such that $\mu_\QQ := \QQ\otimes_\ZZ \mu$ is bijective. In
particular, the element
\[\chi^\rel_{{\rm coh}}(C^*, \mu_\QQ) := [H^0(C^*)] - [\ker(\mu)] +
[H^2(C^*)] - [H^3(C^*)_\codiv] + [\coker(\mu)]  + [H^4(C^*)]\]
belongs to the Grothendieck group $G_0T(\ZZ[G])$ of all finite
$\ZZ[G]$-modules and is a preimage of $\chi_{{\rm coh}}(C^*)$
under the canonical map $G_0T(\ZZ[G]) \ra G_0(\ZZ[G])$ (cf.\
Definition~3.5 in \cite{CKPS}).

The main result of this paper, Theorem~(3.9), implies that for any
nearly perfect complex $(C^*, (L_i)_i, (\tau_i)_i)$ and any $\QQ
[G]$-equivariant `trivialization' isomorphism
\[\lambda: \QQ \otimes H^-(C^*)_\codiv \oplus
\Hom(L_+, \QQ) \,\, \tilde{\ra} \,\, \QQ \otimes H^+(C^*)_\codiv
\oplus \Hom(L_-, \QQ)\] one can define a canonical {\em Euler
characteristic} element $\chi^\rel(C^*,\lambda)$ in the
Grothendieck group $K_0T(\ZZ[G])$ of all finite $\ZZ[G]$-modules
of finite projective dimension which satisfies the following
 key properties: in the obvious commutative diagram
\[\xymatrix{K_0T(\ZZ[G]) \ar[r] \ar[d] & K_0(\ZZ[G]) \ar[d] \\
G_0T(\ZZ[G]) \ar[r] & G_0(\ZZ[G]),}\] $\chi^\rel(C^*, \lambda)$ is
a preimage of $\chi(C^*)$ under the upper horizontal map and, in
the case that $\lambda$ is equal to the map $\mu_\QQ$ described
above, $\chi^\rel(C^*, \lambda)$ is also a preimage of
$\chi^\rel_{{\rm coh}}(C^*,\lambda)$ under the left vertical map.

To define the element $\chi^\rel(C^*,\lambda)$ we proceed as
follows. We first replace $C^*$ by a bounded complex $P^*$ of
projective $\ZZ[G]$-modules which is quasi-isomorphic to $C^*$
(see Lemma~(1.1)). Then for any prime $l$, we show that the
$l$-adic completion $\hat{P}^*$ of $P^*$ is a bounded complex of
cohomologically trivial $\ZZ_l[G]$-modules for which all
cohomology modules are finitely generated. (This follows from
Lemma~(2.1) and Proposition~(2.2).) By a standard procedure (again
see Lemma~(1.1)), we can then find a perfect complex $\tilde{P}^*$
of $\ZZ_l[G]$-modules which is quasi-isomorphic to $\tilde{P}^*$.
Furthermore, the given isomorphism $\lambda$ induces an
isomorphism
\[\lambda_l: \Gr(H^-(\tilde{P}^*))_{\QQ_l} \,\,
\tilde{\ra} \,\, \Gr(H^+(\tilde{P}^*))_{\QQ_l}\] between the
graded objects associated with a certain natural 2-step filtration
on $H^-(\tilde{P}^*)_{\QQ_l}$ and $H^+(\tilde{P}^*)_{\QQ_l}$,
respectively (which we obtain from Proposition~(2.2)). In \S3 we
describe a refined version of the classical construction of
Reidemeister-Whitehead torsion, and we then use this to define a
canonical element $\chi^\rel(\tilde{P}^*,\lambda_l)$ of the
Grothendieck group $K_0T(\ZZ_l[G])$ of all finite
$\ZZ_l[G]$-modules of finite projective dimension. Finally, we
define $\chi^\rel(C^*,\lambda)$ to be equal to the image of the
tuple $(\chi^\rel(\tilde{P}^*,\lambda_l))_{l \, {\rm prime}}$
under the canonical decomposition isomorphism
\[\oplusm_{l \, {\rm prime}} K_0T(\ZZ_l[G]) \cong K_0T(\ZZ[G]).\]

We remark that in several natural arithmetical contexts the
element $\chi^\rel(C^*,\lambda)$ can be directly related
 to the leading terms of associated $L$-functions. For example, in the context
described at the beginning of this introduction, if $k$ is of
characteristic $p$, then in \cite{BuZet} the first named author
has shown that the image of $\chi^\rel(C^*, \mu_\QQ)$ in $K_0T(\ZZ
[G])/K_0T(\ZZ_p[G])$ can be explicitly computed in terms of the
leading term of the $G$-equivariant $L$-function of $X$ at $s =1$.
(The result in \cite{BuZet} is actually for a quasi-projective
variety $X$ of arbitrary dimension). With a view to obtaining
analogous results for varieties over number fields, we have
therefore phrased our definition of $\chi^{\rel}(C^*, \lambda)$ in
\S3 in terms of isomorphisms $\lambda$ which are defined over an
arbitrary field of characteristic zero.

The proof of our main result, Theorem~(3.9), relies upon a certain
  mapping cone construction (explained in \S 1) and
 this in fact yields a natural definition of the Euler characteristic
$\chi(C^*)$ which avoids the somewhat unsatisfactory inductive
procedure used in \cite{CKPS} (see Theorem~(1.3)). We remark that
the idea of using such a mapping cone construction in this context
first occurs in \cite{BuNea}. However, in \cite{BuNea}, only those
 nearly perfect complexes which satisfy a certain natural condition
on homology are considered. Under this condition (which is
  satisfied by all nearly perfect complexes which are known to arise
in arithmetic), it is shown in \cite{BuNea} that the mapping cone
construction which is described in \S1 yields an element of the
derived category which depends (to within isomorphism) only upon
the given nearly perfect complex, and in addition criteria are
described which ensure additivity of the Euler characteristic
$\chi(C^*)$ on distinguished triangles of nearly perfect
complexes.

In this context we remark that it is also possible to shorten some
of the proofs given here by using the language of derived
categories. However, we have decided to eschew such formalism in
order to make this paper accessible to as wide an audience as
possible.

\bigskip

{\em Notations.} Throughout this paper, we fix a finite group $G$,
and all rings are assumed to be unital. The group ring of $G$ over
a ring $R$ is denoted by $R[G]$. As usual, $G_0(R)$ and $K_0(R)$
denote the Grothendieck group of all finitely generated
$R$-modules and the Grothendieck group of all projective finitely
generated $R$-modules, respectively. For any abelian group $A$,
the maximal divisible subgroup of $A$ is denoted by $A_{\rm div}$
and the codivisible quotient $A/A_{\rm div}$ of $A$ by $A_\codiv$.
For any $n\in \NN$, the subgroup of $n$-torsion elements in $A$ is
denoted by ${}_nA$. For a fixed prime $l$, we set $\hat{A}:=
\lim_{\leftarrow_n} A/l^nA$ (the {\em $l$-adic completion of $A$})
and $T_l(A) := \lim_{\leftarrow_n} {}_{l^n}A$ (the {\em $l$-adic
Tate module of $A$}). By a complex $M^*$ in an abelian category we
mean a cochain complex, and the differential from $M^i$ to
$M^{i+1}$ is denoted by $d^i := d^i(M^*)$. In each degree $i$ the
 modules of coboundaries, cocyles and cohomology of $M^*$ are
 denoted by $B^i(M^*)$, $Z^i(M^*)$ and $H^i(M^*)$ respectively. By a
 {\em perfect complex} of $R$-modules we shall mean a bounded complex of finitely
generated projective $R$-modules. For any $\ZZ$-graded object
$M^*$, resp.~$N_*$, we shall write $M^+$ and $M^-$, resp.~$N_+$
and $N_-$, for the direct sum of the objects $M^i$, resp.~$N_i$,
as $i$ runs over all even and odd integers respectively.

\bigskip

\bigskip

\section*{\S 1 The Cone Construction}

In the paper \cite{CKPS}, Chinburg, Kolster, Pappas and Snaith
introduced the notion of a nearly perfect complex, and they
defined its Euler characteristic in $K_0(\ZZ[G])$ by using a
certain somewhat unsatisfactory inductive procedure. The object of
this section is to describe this Euler characteristic as the Euler
characteristic of a certain mapping cone which turns out to be a
bounded complex of cohomologically trivial $\ZZ[G]$-modules whose
cohomology groups are finitely generated. In the first two lemmas,
we will therefore recall how to associate an Euler characteristic
in $K_0(\ZZ[G])$ to any such complex.

{\bf (1.1) Lemma.} {\em Let $R$ be a Dedekind domain and let $C^*$
be a bounded complex of cohomologically trivial $R[G]$-modules.
Then there exists a {\em bounded} complex $P^*$ of projective
$R[G]$-modules
 and a quasi-isomorphism $P^* \ra C^*$. If,
moreover, all cohomology groups $H^i(C^*)$, $i \in \ZZ$, are
finitely generated, then $P^*$ can be chosen to be a perfect
complex of $R[G]$-modules.}

{\em Proof.} It is well-known that there is a complex $Q^*$ of
projective $R[G]$-modules which is bounded to the right together
with a quasi-isomorphism $Q^* \ra C^*$. Let $m \in \NN$ with $C^i
= 0$ for $i <m$. Then the induced homomorphism
\[\tau_{\ge m}(Q^*) := ( \cdots \ra 0 \ra Q^m/d(Q^{m-1}) \ra
Q^{m+1} \ra \cdots) \longrightarrow C^*\] is a quasi-isomorphism,
too. Hence, $Q^m/d(Q^{m-1})$ is cohomologically trivial since the
corresponding mapping cone is acyclic and the modules $C^i$ for
$i\in \ZZ$ and $Q^j$ for $j \ge m+1$, are cohomologically trivial.
So, by Proposition~4.1(b) on p.~457 in \cite{Ch}, there is a short
exact sequence $0\ra P^{m-1} \ra P^m \ra Q^m/d(Q^{m-1}) \ra 0$
where $P^{m-1}$ and $P^m$ are projective $R[G]$-modules. Now the
complex
\[P^* := ( \cdots \ra 0 \ra P^{m-1} \ra P^m \ra Q^{m+1} \ra \cdots
) \] where $P^m \ra Q^{m+1}$ is the composition $\xymatrix{P^m
\ar@{>>}[r] & Q^m/d(Q^{m-1}) \ar[r] & Q^{m+1}}$ is as required.
The second assertion of the lemma is equivalent to asserting that,
under the stated condition on cohomology, each projective module
$P^i$ can be chosen to be finitely generated, and this follows by
a standard argument (see, for example, the proof of Theorem~1.1 on
p.~447 in \cite{Ch} or Proposition 2.2.4 on p.~21 in \cite{Sn}).
\hfill $\square$


{\bf (1.2) Lemma.} {\em Let
\[\xymatrix{ P^* \ar[r]^{\alpha^*}& C^* \ar[d]^{f^*} \\
Q^* \ar[r]^{\beta^*}& D^*}\] be a diagram of complexes in an
abelian category $\cM$ where $P^*$ is a complex of projective
objects which is bounded to the right and where $\beta^*$ is a
quasi-isomorphism. Then there exists a homomorphism of complexes
$h^*: P^* \ra Q^*$ such that $f^*\alpha^*$ is homotopic to
$\beta^* h^*$. In particular, in the situation of the second
statement of Lemma~(1.1), the element
\[\chi(C^*) := \sum_{i \in \ZZ} (-1)^i [P^i] \in K_0(R[G])\]
does not depend on the chosen quasi-isomorphism $P^* \ra C^*$, and
we call it the {\em Euler characteristic of $C^*$}.}

{\em Proof.} The first statement is a standard result (see, for
example, the proof of Proposition 2.2.4(ii) on p.~21 in
\cite{Sn}). If, in the situation of the second statement of
Lemma~(1.1), we choose another perfect complex $Q^*$ together with
a quasi-iso\-morphism $Q^* \ra C^*$, then we obtain a
quasi-isomorphism $P^* \ra Q^*$ by the first statement of
Lemma~(2.1). The corresponding mapping cone being acyclic, we
obtain
\[\sum_{i\in \ZZ} (-1)^i [P^i] - \sum_{i \in \ZZ} (-1)^i [Q^i] = \sum_{i
\in \ZZ} (-1)^i [P^i \oplus Q^{i-1}] = 0 \quad \textrm{in} \quad
K_0(R[G]),\] as was to be shown. \hfill $\square$

Let now $(C^*, (L_i)_i, (\tau_i)_i)$ be a fixed nearly perfect
complex: thus, $C^*$ is a bounded complex of cohomologically
trivial $\ZZ[G]$-modules and, for each~$i\in \ZZ$, the
$\ZZ[G]$-module $H^i(C^*)_\codiv$ is finitely generated , $L_i$ is
a $\ZZ$-free, finitely generated $\ZZ[G]$-module, and $\tau_i:
\Hom_\ZZ(L_i, \QQ/\ZZ) \ra H^i(C^*)_{\rm div}$ is a
$\ZZ[G]$-equivariant isomorphism.

We recall that a uniquely divisible $\ZZ[G]$-module is
cohomologically trivial and therefore (as observed earlier) admits
a $\ZZ[G]$-projective resolution of length at most $1$. For
 each $i$, we fix such a resolution
\[ \xymatrix{0 \ar[r] & R^{i-1} \ar@{^{(}->}[r] & Q^i
\ar[r]^>>>>{\varepsilon^i}  & \Hom(L_i, \QQ) \ar[r] & 0 } \] of
$\Hom(L_i,\QQ)$. Furthermore, we choose a map of complexes
\[\xymatrix{ \cdots \ar[r]& 0 \ar[r] \ar[d] & R^{i-1} \ar[r]
\ar[d]^{\beta^{i-1}} & Q^i \ar[r] \ar[d]^{\alpha^i} & 0 \ar[r]
\ar[d] & \cdots \\ \cdots \ar[r] & C^{i-2} \ar[r] & C^{i-1} \ar[r]
& C^i \ar[r] & C^{i+1} \ar[r] & \cdots }\] such that the following
diagrams commute:
\[\xymatrix{Q^i \ar@{>>}[r]^>>>>{\varepsilon^i} \ar@{-->}[ddd]^{\alpha^i}
& \Hom(L_i, \QQ) \ar@{>>}[d] \\ & \Hom(L_i,\QQ/\ZZ)
\ar[d]^{\tau_i}_\wr \\ & H^i(C^*)_{\rm div} \ar@{^{(}->}[d] \\
Z^i(C^*) \ar@{>>}[r]  & H^i(C^*) } \quad \textrm{ and } \quad
\xymatrix{R^{i-1} \ar@{^{(}->}[r] \ar[ddr]
\ar@{-->}[dd]^{\beta^{i-1}} & Q^i \ar[ddr]^{\alpha^i} \\ \\C^{i-1}
\ar@{>>}[r] & B^i(C^*) \ar@{^{(}->}[r] & Z^i(C^*). } \] These maps
induce a homomorphism of complexes
\[ \xymatrix{ \cdots \ar[r] & Q^{i-1} \oplus R^{i-1} \ar[r]
\ar[d]^{(\alpha^{i-1}, \beta^{i-1})} & Q^i \oplus R^i \ar[r]
\ar[d]^{(\alpha^i,\beta^i)} & Q^{i+1}\oplus R^{i+1} \ar[r]
\ar[d]^{(\alpha^{i+1},\beta^{i+1})} & \cdots \\ \cdots \ar[r]&
C^{i-1} \ar[r] & C^i \ar[r] & C^{i+1} \ar[r] & \cdots } \] where
the upper complex (in the sequel denoted by $(Q\oplus R)^*$) is
the direct sum of the complexes $ \cdots \ra 0 \ra R^{i-1} \ra Q^i
\ra 0 \ra \cdots$, $i \in \ZZ$. Let $\Cone_C^*$ denote the mapping
cone of this homomorphism. Then we have a short exact sequence of
complexes
\[ 0 \ra C^* \ra \Cone_C^* \ra (Q\oplus R)^*[1] \ra 0\]
which yields the long exact sequence
\[ \xymatrix{\cdots \ar[r] & \Hom(L_i,\QQ) \ar[r]^{\partial} &
H^i(C^*) \ar[r] & H^i(\Cone_C^*) \ar[r]& \Hom(L_{i+1}, \QQ) \ar[r]
& \cdots }\] where the connecting homomorphism $\partial$ is the
composition \[\xymatrix{\Hom(L_i, \QQ) \ar@{>>}[r] &
\Hom(L_i,\QQ/\ZZ) \ar[r]^<<<{\tau_i}_<<<\sim & H^i(C^*)_{\rm div}
\ar@{^{(}->}[r] & H^i(C^*) }.\] Therefore we obtain natural short
exact sequences
\begin{equation}
0 \ra H^i(C^*)_\codiv \ra H^i(\Cone_C^*) \ra \Hom(L_{i+1},\ZZ) \ra
0, \quad i \in \ZZ.
\end{equation}
In particular, the $\ZZ[G]$-module $H^i(\Cone_C^*)$ is finitely
generated for all $i$. Since all of the $\ZZ[G]$-modules
$\Cone_C^i$, $i \in \ZZ$, are cohomologically trivial, we obtain
an Euler characteristic $\chi(\Cone_C^*) \in K_0(\ZZ[G])$ by using
Lemma~(1.2).

{\bf (1.3) Theorem.} {\em The class $\chi(\Cone_C^*) \in
K_0(\ZZ[G])$ is equal to the class $\chi(C^*):= \chi((C^*,
(L_i)_i, (\tau_i)_i))$ which is defined in (\cite{CKPS},
Definition~2.11).}

{\em Proof.} First, we assume that at most one module $C:= C^n$ in
the complex $C^*$ is non zero. Then $\chi(C^*)$ is defined in the
following way. We choose a finitely generated projective
$\ZZ[G]$-module $F$ together with a $\ZZ[G]$-module homomorphism
$F\ra C$ such that the composition $F \ra C \ra C_\codiv$ is
surjective. This homomorphism yields the following commutative
diagram with exact rows and columns (where $L:= L_n$ and $M$ and
$N$ are defined in such a way that the two right columns become
exact):
\[\xymatrix{&0 \ar[d] & 0 \ar[d] & 0 \ar[d] \\
0 \ar[r] & \Hom(L,\ZZ) \ar[r] \ar[d] & M \ar[r] \ar[d] & N \ar[r]
\ar[d] & 0 \\ 0 \ar[r] & \Hom(L,\QQ) \ar[r] \ar[d] & \Hom(L,\QQ)
\oplus F \ar[r] \ar[d] & F\ar[r] \ar[d] & 0 \\ 0 \ar[r] &
\Hom(L,\QQ/\ZZ) \ar[r] \ar[d] & C \ar[r] \ar[d] & C_\codiv \ar[r]
\ar[d]& 0 \\ & 0 & 0 & 0 } \] Then $M$ is a finitely generated
projective $\ZZ[G]$-module (see Corollary~2.8 in \cite{CKPS}) and,
by definition, $\chi(C^*) = (-1)^n([F]-[M])$ in $K_0(\ZZ[G])$ or,
in other words, $\chi(C^*)$ is the Euler characteristic of the
perfect complex $\cdots \ra 0 \ra M \ra F \ra 0 \ra \cdots$ (with
$F$ in degree $n$, of course) resulting from the diagram above.
This diagram also shows that this complex is quasi-isomorphic to
the complex $\cdots \ra 0 \ra \Hom(L,\QQ) \ra C \ra 0 \ra \cdots$.
Since we furthermore have an obvious quasi-isomorphism from the
complex $\Cone_C^*$ to the complex $\cdots \ra 0 \ra \Hom(L,\QQ)
\ra C \ra 0 \ra \cdots$, we obtain a quasi-isomorphism from the
complex $\cdots \ra 0 \ra M \ra N \ra 0 \ra \cdots$ to the complex
$\Cone_C^*$ (by Lemma (1.2)) which
proves that $\chi(C^*) = \chi(\Cone_C^*)$, as was to be shown. \\
We now proceed by induction on the length of $C^*$. Let $n\in \ZZ$
such that $C^n \not= 0$ and $C^r = 0$ for all $r>n$. We choose a
finitely generated projective $\ZZ[G]$-module $F^n$ together with
a homomorphism $F^n \ra C^n$ such that the composition
\[\xymatrix{F^n \ar[r] &C^n \ar@{>>}[r] &H^n(C^*) \ar@{>>}[r]&
H^n(C^*)_\codiv}\]  is surjective. Furthermore, we choose a
projective $\ZZ[G]$-module $S^n$ together with an epimorphism
$\xymatrix{S^n \ar@{>>}[r] & B^n(C^*)}$ Then we obtain the
following commutative diagram
\[\xymatrix{P^n := S^n \oplus Q^n \oplus F^n
\ar@{>>}[rr]^{0 \oplus \varepsilon^n \oplus \id} \ar@{>>}[d] &&
\Hom(L^n,\QQ)\oplus F^n \ar@{>>}[d] \\ C^n \ar@{>>}[rr] &&
H^n(C^*). }\] Let the $\ZZ[G]$-modules $B_P^n$ and $M^n$ be
defined by the following commutative diagram with exact rows and
columns:
\[\xymatrix{& 0 & 0 \\
0 \ar[r] & M^n \ar[u] \ar[r] & \Hom(L^n,\QQ) \oplus F^n \ar[u]
\ar@{>>}[r]& H^n(C^*) \ar[r] & 0\\ 0 \ar[r] & B_P^n \ar[u] \ar[r]
& P^n \ar@{>>}[u] \ar@{>>}[r] & H^n(C^*) \ar@{=}[u] \ar[r] & 0 \\
& S^n \oplus R^{n-1} \ar[u] \ar@{=}[r]& S^n \oplus R^{n-1} \ar[u]
\\ & 0 \ar[u] & 0. \ar[u] }\]  Let $E^{n-1}$ be defined by the
following pull-back diagram where the epimorphism $\xymatrix{B_P^n
\ar@{>>}[r] & B^n(C^*)}$ is induced by the epimorphism
$\xymatrix{P^n = S^n \oplus Q^n \oplus F^n \ar@{>>}[r] & C^n}$
introduced above:
\[\xymatrix{0 \ar[r] &Z^{n-1}(C^*) \ar[r] \ar@{=}[d] & E^{n-1} \ar[r]
\ar@{>>}[d] & B_P^n \ar[r] \ar@{>>}[d] & 0 \\ 0 \ar[r]&
Z^{n-1}(C^*) \ar[r] & C^{n-1} \ar[r]& B^n(C^*) \ar[r]& 0. } \]
Since $S^n \oplus R^{n-1}$ is $\ZZ[G]$-projective, we can lift
the inclusion $S^n \oplus R^{n-1} \hookrightarrow B_P^n$ to an
inclusion $S^n \oplus R^{n-1} \hookrightarrow E^{n-1}$ and define
$D^{n-1} := E^{n-1}/(S^n \oplus R^{n-1})$. Then we obtain the
following exact sequence:
\[0 \ra Z^{n-1}(C^*) \ra D^{n-1} \ra M^n \ra 0.\]
Putting all diagrams together, we obtain the following commutative
diagram with exact rows and vertical epimorphisms:
\begin{equation}
\xymatrix{0 \ar[r] & Z^{n-1}(C^*) \ar[r] & D^{n-1} \ar[r] &
\Hom(L_n,\QQ) \oplus F^n \ar[r] & H^n(C^*) \ar[r]& 0 \\ 0 \ar[r]&
Z^{n-1}(C^*) \ar@{=}[u] \ar[r] \ar@{=}[d] & E^{n-1} \ar@{>>}[u]
\ar[r] \ar@{>>}[d] & P^n \ar@{>>}[u] \ar[r] \ar@{>>}[d] & H^n(C^*)
\ar@{=}[u] \ar[r] \ar@{=}[d] & 0 \\ 0 \ar[r] & Z^{n-1}(C^*) \ar[r]
& C^{n-1} \ar[r] & C^n \ar[r] & H^n(C^*) \ar[r] & 0. }
\end{equation}
Since the kernel of the epimorphism $\xymatrix{E^{n-1} \ar@{>>}[r]
& C^{n-1}}$ is isomorphic to the kernel of the epimorphism
$\xymatrix{P^n \ar@{>>}[r] & C^n}$, $E^{n-1}$ is cohomologically
trivial; hence $D^{n-1}$ is also cohomologically trivial. The
diagram~(2) furthermore shows that the 1-extension
\[0 \ra Z^{n-1}(C^*) \ra D^{n-1} \ra M^n \ra 0\] is a preimage (which
is unique by the proof of Lemma~2.4 in \cite{CKPS}) of the
tautological 2-extension
\[ 0 \ra Z^{n-1}(C^*) \ra C^{n-1} \ra C^n \ra H^n(C^*) \ra 0\]
under the connecting homomorphism
\[\Ext^1_{\ZZ[G]}(M^n, Z^{n-1}(C^*)) \ra
\Ext^2_{\ZZ[G]}(H^n(C^*),Z^{n-1}(C^*))\] which is associated with
the short exact sequence
\[0 \ra M^n \ra \Hom(L^n,\QQ) \oplus F^n \ra H^n(C^*) \ra
0.\] Let $D^*$ denote the complex
\[\cdots \ra C^{n-3} \ra C^{n-2} \ra D^{n-1} \ra 0 \ra \cdots\]
where $C^{n-2} \ra D^{n-1}$ is the composition $C^{n-2} \ra
Z^{n-1}(C^*) \ra D^{n-1}$. Then, by definition, we have
\[\chi(C^*) =
\chi(D^*) + (-1)^n[F^n] \quad \textrm{in} \quad K_0(\ZZ[G])\]
where we view $D^*$ as a nearly perfect complex (of a smaller
length than $C^*$) as in Corollary~2.10 in \cite{CKPS}. Let
$D^*_{\rm aug}$ respectively $E^*$ denote the complexes
\[ \cdots \ra C^{n-3} \ra C^{n-2} \ra D^{n-1} \ra \Hom(L_n,\QQ) \oplus F^n \ra
0 \ra \cdots\] respectively
\[\cdots \ra C^{n-3} \ra C^{n-2} \ra E^{n-1} \ra P^n \ra 0 \ra \cdots \]
resulting from the diagram (2). From the construction we obtain
natural homomorphisms of complexes
\[(Q\oplus R)^* \ra E^* \quad \textrm{and} \quad
(Q\oplus R)^* \ra D^*_{\rm aug} \]  such that the following
triangles commute:
\[\xymatrix{&(Q\oplus R)^* \ar[dl] \ar[dr]^{(\alpha^i,\beta^i)_i}
\\ E^* \ar[rr] && C^* } \quad \textrm{and} \quad \xymatrix{&(Q\oplus
R)^* \ar[dl] \ar[dr] \\ E^* \ar[rr] && D^*_{\rm aug} } \]
Therefore, we obtain quasi-isomorphisms between the corresponding
cones:
\[\Cone^*_C \leftarrow \Cone^*_E \ra \Cone^*_{D_{\rm aug}}\]
hence:
\[ \chi(\Cone_C^*) = \chi(\Cone^*_{D_{\rm aug}}) \quad \textrm{in} \quad
K_0(\ZZ[G]).\] Furthermore, the natural epimorphism from
$\Cone^*_{D_{\rm aug}}$ to the mapping cone $M^*$ of the
homomorphism of complexes
\[\xymatrix{\cdots \ar[r] & Q^{n-2} \oplus R^{n-2} \ar[r] \ar[d] &
Q^{n-1} \ar[r] \ar[d] & 0 \ar[r] \ar[d] & \cdots \\ \cdots \ar[r]&
D^{n-2} \ar[r]& D^{n-1} \ar[r] & F^n \ar[r] & \cdots }\]  is a
quasi-isomorphism since the kernel
\[\cdots \ra 0 \ra R^{n-1} \ra Q^n \ra \Hom(L_n,\QQ) \ra 0 \ra
\cdots\] is acyclic. Finally, by using the inductive hypothesis it
is easy to see that the Euler characteristic of $M^*$ is equal to
$\chi(D^*) + (-1)^n[F^n]$. So, Theorem (1.3) is proved:
\[\chi(\Cone_C^*) = \chi(\Cone_{D_{\rm aug}}^*) = \chi(M^*) =
\chi(D^*) + (-1)^n[F^n] = \chi(C^*) \quad \textrm{in} \quad
K_0(\ZZ[G]).\] \hspace*{\fill} $\square$

\bigskip

\section*{\S 2 The $l$-adic Euler Characteristic of a Nearly\\
Perfect Complex}

In this section we fix a prime $l$. We will show that the image of
the Euler characteristic $\chi(C^*)$ of a nearly perfect complex
$C^*$ in the Grothendieck group $K_0(\ZZ_l[G])$ of all finitely
generated projective $\ZZ_l[G]$-modules is equal to the Euler
characteristic of any complex which is obtained by replacing $C^*$
with a quasi-isomorphic complex of projective $\ZZ[G]$-modules and
then passing to $l$-adic completions.

We begin with the following easy observation.

{\bf (2.1) Lemma.} {\em Let $P$ be a $\ZZ$-torsion-free
cohomologically trivial $\ZZ[G]$-module. Then its $l$-adic
completion $\hat{P}$ is also cohomologically trivial.}

{\em Proof.} Since the transition maps in the inverse system
$(P/l^nP)_{n\ge 0}$ are surjective, we have a short exact sequence
\[0 \ra \limm_{\stackrel{\longleftarrow}{{n}}} P/l^nP \ra \prod_n P/l^nP \ra \prod_n
P/l^nP \ra 0\] where the right homomorphism maps a tuple
$(x_n)_{n\in \NN}$ to the tuple $(x_n - \overline{x_{n+1}})_{n \in
\NN}$. Now, let $\hat{H}$ denote Tate cohomology with respect to
some subgroup of $G$. Then we obtain a long exact sequence
\[\cdots \ra \hat{H}^s(\hat{P}) \ra \prod_n \hat{H}^s(P/l^nP) \ra
\prod_n \hat{H}^s(P/l^nP) \ra \hat{H}^{s+1} (\hat{P}) \ra \cdots
\]
Since $P$ is $\ZZ$-torsion-free and cohomologically trivial,
$\hat{H}^s(P/l^nP)$ vanishes for all $s \in \ZZ$ and $n\in \NN$.
Thus, $\hat{H}^s(\hat{P}) =0$ for all $s\in \ZZ$, as was to be
shown. \hfill $\square$

{\bf (2.2) Proposition.} {\em Let $P^*$ be a complex of
$\ZZ$-torsion-free $\ZZ[G]$-modules such that the $l^n$-torsion
subgroup ${}_{l^n}H^i(P^*)$ and the $l^n$-quotient group $H^i(P^*)
\otimes \ZZ/l^n\ZZ$ are finite for all $i\in \ZZ$ and all $n \in
\NN$. Then there is a natural short exact sequence
\[0 \ra \widehat{H^i(P^*)_\codiv} \ra H^i(\hat{P}^*) \ra
T_l(H^{i+1}(P^*)_{{\rm div}}) \ra 0 \] of $\ZZ_l[G]$-modules for
all $i\in \ZZ$.}

{\em Proof.} The exact sequences
\[0 \ra P^* \ra P^* \ra P^*\otimes \ZZ/l^n\ZZ \ra 0, \quad n\ge
0,\] yield the familiar short exact sequences
\[0 \ra H^i(P^*) \otimes \ZZ/l^n\ZZ \ra H^i(P^* \otimes
\ZZ/l^n\ZZ) \ra {}_{l^n}H^{i+1}(P^*) \ra 0, \quad n\in \NN, \quad
i\in \ZZ.\] Since these are short exact sequences of finite
groups, they remain exact upon completion by the Mittag-Leffler
criterion. Thus we have short exact sequences
\begin{equation}0\ra \widehat{H^i(P^*)} \ra \limm_{\stackrel{\longleftarrow}{n}}
H^i(P^* \otimes \ZZ/l^n\ZZ) \ra T_l(H^{i+1} (P^*)) \ra 0, \quad
i\in \ZZ.
\end{equation}
The exact sequence
\[0 \ra H^i(P^*)_{{\rm div}} \ra H^i(P^*) \ra H^i(P^*)_\codiv \ra 0\]
yields an isomorphism $H^i(P^*) \otimes \ZZ/l^n\ZZ \cong
H^i(P^*)_\codiv \otimes \ZZ/l^n\ZZ$ for all $n$, hence
\begin{equation}
\widehat{H^i(P^*)} \cong \widehat{H^i(P^*)_\codiv} \textrm{ for
all } i \in \ZZ;
\end{equation}
furthermore, we obtain left exact sequences
\[0 \ra {}_{l^n}(H^i(P^*)_{\rm div}) \ra {}_{l^n}H^i(P^*) \ra
{}_{l^n}(H^i(P^*)_\codiv), \quad i\in \ZZ, \quad n \in \NN;\]
hence
\begin{equation}
T_l(H^i(P^*)_{\rm div}) \cong T_l(H^i(P^*)) \textrm{ for all } i
\in \ZZ
\end{equation}
since, for any abelian group $A$, $T_l(A_\codiv)$ vanishes.
(Proof: if $T_l(A_\codiv)$ contained a non-zero element
$(x_n)_{n\in \NN}$, then the subgroup $U$ of $A_\codiv$ generated
by the $l^\infty$-torsion elements $x_n$, $n\in \NN$, would be
non-zero
 and divisible (by {\em all} primes), hence the preimage of $U$
under the canonical epimorphism $A\ra A_\codiv$ would be a
divisible subgroup of $A$ bigger than $A_{\rm div}$.)  \\
Furthermore, we have a short exact sequence of complexes
\[0 \ra \hat{P}^* \ra \prod_n P^* \otimes \ZZ/l^n\ZZ \ra \prod_n
P^* \otimes \ZZ/l^n\ZZ \ra 0\] where the right map is defined as
in the proof of Lemma~(2.1). Thus we obtain the long exact
sequence
\[\cdots \ra H^i(\hat{P}^*) \ra \prod_n H^i(P^*\otimes \ZZ/l^n\ZZ)
\ra \prod_n H^i(P^* \otimes \ZZ/l^n\ZZ) \ra H^{i+1} (\hat{P}^*)
\ra \cdots \] which yields the exact Milnor sequence
\[0 \ra {\limm_{\stackrel{\longleftarrow}{n}}}^1 H^{i-1}(P^* \otimes \ZZ/l^n\ZZ) \ra
H^i(\hat{P}^*) \ra \limm_{\stackrel{\longleftarrow}{n}} H^i(P^*
\otimes \ZZ/l^n\ZZ) \ra 0.\] Since $H^{i-1}(P^*\otimes
\ZZ/l^n\ZZ)$ is finite for all $n$, the inverse system
$(H^{i-1}(P^* \otimes \ZZ/l^n\ZZ))_{n\ge 0}$ satisfies the
Mittag-Leffler condition, so
${\limm_{\stackrel{\longleftarrow}{n}}}^1 H^{i-1}(P^* \otimes
\ZZ/l^n\ZZ) = 0$, and hence
\begin{equation}
H^i(\hat{P}^*) \cong \limm_{\stackrel{\longleftarrow}{n}} H^i(P^*
\otimes \ZZ/l^n\ZZ) \textrm{ for all } i.
\end{equation}
Upon combining (3), (4), (5) and (6) we obtain Proposition~(2.2).
\hfill $\square$

{\bf (2.3) Corollary.} {\em Let $P^* \ra Q^*$ be a
quasi-isomorphism of complexes as in Proposition~(2.2). Then the
induced morphism $\hat{P}^* \ra \hat{Q}^*$ is also a
quasi-isomorphism.}

{\em Proof.} Obvious. \hfill $\square$

Now, let $C^* = (C^*, (L_i)_i, (\tau_i)_i)$ be a nearly perfect
complex of $\ZZ[G]$-modules. We choose a bounded complex $P^*$ of
projective $\ZZ[G]$-modules together with a quasi-isomorphism $P^*
\ra C^*$ as in Lemma~(1.1). Then the $l$-adic completion
$\hat{P}^*$ is a complex of cohomologically trivial
$\ZZ_l[G]$-modules by Lemma~(2.1). By definition, we have short
exact sequences
\[0 \ra \Hom_\ZZ(L^i, \QQ/\ZZ) \ra H^i(P^*) \ra H^i(C^*)_\codiv
\ra 0, \quad i \in \ZZ, \] where $L_i$ is a $\ZZ$-free finitely
generated $\ZZ[G]$-module and the $\ZZ[G]$-module
$H^i(C^*)_\codiv$ is finitely generated. So, the conditions of
Proposition~(2.2) are satisfied and we obtain the short exact
sequences
\[0 \ra H^i(C^*)_\codiv \otimes \ZZ_l \ra H^i(\hat{P}^*) \ra
\Hom(L^{i+1}, \ZZ_l)\ra 0, \quad i\in \ZZ,\] since $T_l(\QQ/\ZZ) =
\ZZ_l$ and since $l$-adic completion of a finitely generated
$\ZZ[G]$-module is the same as tensoring with $\ZZ_l$. In
particular, the $\ZZ_l[G]$-modules $H^i(\hat{P}^*)$, $i\in \ZZ$,
are finitely generated. Thus, by Lemma~(1.1), there is a perfect
complex $\tilde{P}^*$ of $\ZZ_l[G]$-modules together with a
quasi-isomorphism $\tilde{P}^* \ra \hat{P}^*$. If we choose
another quasi-isomorphism $Q^* \ra C^*$ as in Lemma~(1.1) and a
corresponding quasi-isomorphism $\tilde{Q}^* \ra \hat{Q}^*$ as
above, then we obtain a quasi-isomorphism $\tilde{P}^* \ra
\tilde{Q}^*$ by applying Lemma~(1.2) twice and by using
Corollary~(2.3). In particular therefore, the class
\[\chi_l(C^*) := \sum_{i \ge 0} (-1)^i [\tilde{P}^i] \in
K_0(\ZZ_l[G])\] does not depend on the above choices. We call
$\chi_l(C^*) \in K_0(\ZZ_l[G])$ the {\em $l$-adic Euler
characteristic of the nearly perfect complex
$(C^*,(L_i)_i,(\tau_i)_i)$}.

{\bf (2.4) Theorem.} {\em Under the canonical homomorphism
$K_0(\ZZ[G]) \ra K_0(\ZZ_l[G])$, the element $\chi(C^*) =
\chi((C^*, (L_i)_i, (\tau_i)_i))$ of $K_0(\ZZ[G])$ (see
Definition~2.6 in \cite{CKPS}) maps to $\chi_l(C^*)$ in
$K_0(\ZZ_l[G])$.}

{\em Proof.} As in \S1, for each $i$, we choose a
$\ZZ[G]$-projective resolution
\[\xymatrix{0 \ar[r]& R^{i-1} \ar@{^{(}->}[r]& Q^i
\ar[r]^>>>>{\varepsilon^i} & \Hom(L_i,\QQ) \ar[r] & 0 } \] and a
map of complexes $(Q\oplus R)^* \ra P^*$ such that the composition
with the quasi-isomorphism $P^* \ra C^*$ chosen above is of the
form considered in \S1. Let $\Cone_P^*$ denote the corresponding
mapping cone. The corresponding short exact sequence of complexes
\[0 \ra P^* \ra \Cone_P^* \ra (Q\oplus R)^*[1] \ra 0\]
yields a quasi-isomorphism between the complexes $P^*/l^n P^*$ and
$\Cone_P^*/l^n \Cone_P^*$ for all $n\ge 0$ since the cohomology of
$(Q\oplus R)^*$ is uniquely divisible. Hence we obtain a
quasi-isomorphism $\hat{P}^* \ra \widehat{\Cone}_P^*$ between the
$l$-adically completed complexes (by means of the isomorphism~(6)
in the proof of Proposition~(2.2)). Furthermore, we choose a
quasi-isomorphism $P_{\Cone}^* \ra \Cone_P^*$ from a perfect
complex $P_\Cone^*$ to $\Cone_P^*$ (which exists by Lemma~(1.1)).
Then the induced homomorphism $\widehat{P_\Cone}^* \ra
\widehat{\Cone}_P^*$ between the completed complexes is a
quasi-isomorphism by Corollary~(2.3). Hence, by Lemma~(1.2), we
obtain a quasi-isomorphism from the perfect complex
$\widehat{P_\Cone}^*$ to the complex $\hat{P}^*$. Therefore,
$\chi_l(C^*) \,\, \stackrel{{\rm def}}{=} \,\, \chi(\hat{P}^*)
\,\, \stackrel{{\rm def}}{=} \,\, \chi(\widehat{P_\Cone}^*)$ is
the image of $\chi(P_\Cone^*)$ under the canonical map
$K_0(\ZZ[G]) \ra K_0(\ZZ_l[G])$. Finally, by using the composition
of the quasi-isomorphism $P_\Cone^* \ra \Cone_P^*$ with the
obvious quasi-isomorphism $\Cone_P^* \ra \Cone_C^*$ we deduce from
Theorem~(1.3) that $\chi(P_\Cone^*)$ is equal to the Euler
characteristic of the nearly perfect complex $C^*$. So,
Theorem~(2.4) is proved. \hfill $\square$

\bigskip

\section*{\S 3 Refined Euler Characteristics}

Let $C^* = (C^*, (L_i)_i, (\tau_i)_i)$ be a nearly perfect
complex. In this section we also assume given a field $E$ of
characteristic $0$ and an $E[G]$-equivariant {\em trivialization}
isomorphism
\begin{equation}\label{Etriv} \lambda: E \otimes H^-(C^*)_\codiv
\oplus \Hom(L_+,E) \,\, \tilde{\ra} \,\, E \otimes H^+(C^*)_\codiv
\oplus \Hom(L_-,E).\end{equation} We shall refer to any such pair
$(C^*,\lambda)$ as an {\em $E$-trivialized nearly perfect
complex}.

The main aim of this section (achieved in Theorem~(3.9)) is to
associate to each such pair
 $(C^*,\lambda)$ a canonical Euler characteristic element
 $\chi^{\rm rel}(C^*,\lambda)$. This element belongs to
 the relative algebraic $K$-group $K_0(\ZZ [G],E)$ (whose definition is recalled below)
and constitutes a natural refinement of the Euler characteristic
$\chi(C^*)$ of $C^*$ which has been discussed in \S1.

Before proceeding we reassure the reader that $E$-trivialized
nearly perfect complexes $(C^*,\lambda)$ arise naturally in many
arithmetical contexts (with the isomorphism $\lambda$ arising
 via regulator maps or height pairings which are defined over
 $E=\RR$ or $E=\QQ_l$). Further, in many such cases it can be shown that the element
 $\chi^{\rel}(C^*,\lambda)$ provides an important means of
 relating natural Euler characteristics to special values of associated
 $L$-functions (for more details in this direction see \cite{BuNea}).

The initial constructions and results of this section are modeled
on those of \S1 in \cite{BuIwa} (as reviewed in \S1.2 of
\cite{BuComp}). However, it will be convenient for us to consider
the following more general context. Let $\varphi: R \ra S$ be a
homomorphism between unital rings. We assume that the tensor
functor $S \otimes_R -$ from the category of $R$-modules to the
category of $S$-modules is exact. (We will simply write $M_S$ for
$S \otimes_R M$ for any $R$-module $M$.) We further assume that
the ring $S$ is semisimple, i.e., that all $S$-modules are
projective (and hence also injective).

We recall that the {\em relative Grothendieck group $G_0(R,
\varphi)$ of coherent $R$-modules} is defined to be the abelian
group with generators $[A,g,B]$ where $A$, $B$ are finitely
generated $R$-modules and $g$ is an $S$-module isomorphism from
$A_S$ to $B_S$ and the following relations:\\ (Ri) $[A,g,B] = [A',
g', B'] + [A'', g'', B'']$ whenever there is a short exact
sequence of triples (with the obvious meaning)
\[0 \ra (A', g', B') \ra (A,g,B) \ra (A'',g'',B'') \ra 0.\]
(Rii) $[A,hg,C]= [A,g,B] + [B,h,C]$.\\ In particular, we have
$[A,\id,A] = 0$ and $[A,g,B] = -[B,g^{-1},A]$ in $G_0(R,\varphi)$.

We now describe the main example which we have in mind.

{\em (3.1) Example.} Let $E$ be any field of characteristic 0, and
let $\varphi:\ZZ[G] \ra E[G]$ be the canonical inclusion of group
rings. In this case, we also write $G_0(\ZZ[G],E)$ for
$G_0(\ZZ[G],\varphi)$. It is well-known and easy to prove that the
association $[M] \mapsto [0,0,M]$ induces a well-defined
isomorphism
\[G_0T(\ZZ[G]) \,\, \tilde{\ra} \,\, G_0(\ZZ[G],\QQ)\]
where $G_0T(\ZZ[G])$ denotes the Grothendieck group of all
$\ZZ$-torsion $\ZZ[G]$-modules; its inverse is given by sending a
generator $[A,g,B]$ as above to the element
\[[\coker(h)] - [{\rm ker}(h)] - [B/nB] + [{}_nB]\]
where $h \in \Hom_{\ZZ[G]}(A,B)$ and $n\in \ZZ$ are chosen in such
a way that $h/n = g$ in $\Hom_{\QQ[G]}(A_\QQ,B_\QQ) =
\Hom_{\ZZ[G]}(A,B)_\QQ$. (Here, we write $A_\QQ$ for $\QQ
\otimes_\ZZ A \cong \QQ[G] \otimes_{\ZZ[G]} A$.) \\ (This
statement can of course be generalized to the situation in which
$\ZZ$ is replaced by any Dedekind domain and $\QQ$ by the
corresponding field of fractions.)

Let now $M^*$ be a bounded complex of finitely generated
$R$-modules. We assume that, for each $i$, we are given a (finite,
exhaustive and separated) decreasing filtration $(F^nH^i)_{i \in
\ZZ}$ on the cohomology modules $H^i := H^i(M^*)$ and that we are
given an isomorphism
\[\lambda: \Gr(H^-)_S \,\, \tilde{\ra} \,\, \Gr(H^+)_S\]
from the module $\Gr(H^-)_S := \oplus_{i \,\, \odd} \oplus_{n \in
\ZZ} \Gr^n(H^i)_S := \oplus_{i \, \odd} \oplus_{n \in \ZZ}
(F^nH^i/F^{n+1}H^i)_S$ to the module $\Gr(H^+)_S := \oplus_{i \,\,
\even} \oplus_{n \in \ZZ} (F^nH^i/F^{n+1}H^i)_S$. By the
assumption on $\varphi$ and $S$, we have split short exact
sequences
\[0 \ra Z^i_S \ra M^i_S \ra B^{i+1}_S \ra 0, \quad
0 \ra B^i_S \ra Z^i_S \ra H^i_S \ra 0\] and
\[0 \ra (F^{n+1}H^i)_S \ra (F^nH^i)_S \ra \Gr^n(H^i)_S \ra 0, \quad n \in
\ZZ,\] where we write $B^i := B^i(M^*)$, $Z^i := Z^i(M^*)$ and
$H^i:= H^i(M^*)$ for brevity. So we obtain an $S$-module
isomorphism
\begin{eqnarray*}
\lefteqn{\lambda_{M^*}: M^-_S := \oplusm_{i \,\, \odd} M^i_S \,\,
\stackrel{(\iota_1)}{\cong} \,\, \oplusm_{i \,\, \odd} (Z^i_S
\oplus B^{i+1}_S)} \\&& \stackrel{(\iota_2)}{\cong} \,\,
\oplusm_{i \,\, \odd} (H^i_S \oplus B^i_S \oplus B^{i+1}_S) \,\,
\stackrel{(\iota_3)}{\cong} \,\, \oplus_{i \, \odd} (\oplus_{n\in
\ZZ} (\Gr^n(H^i)_S) \oplus B^i_S \oplus B^{i+1}_S)  \\ &&
\stackrel{(\iota_4)}{\cong} \,\, \oplus_{i \, \even} (\oplus_{n
\in \ZZ} (\Gr^n(H^i)_S) \oplus B^i_S \oplus B^{i+1}_S) \,\,
\stackrel{(\iota_5)}{\cong} \,\, \oplus_{i \,\, \even}(H^i_S
\oplus B^i_S \oplus B^{i+1}_S)
\\ && \stackrel{(\iota_6)}{\cong} \,\, \oplusm_{i \,\, \even}(Z^i_S
\oplus B^{i+1}_S) \,\, \stackrel{(\iota_7)}{\cong} \,\, \oplusm_{i
\,\, \even} M^i_S =: M^+_S
\end{eqnarray*}
where the isomorphisms ($\iota_1$), ($\iota_2$), ($\iota_3$),
($\iota_5$), ($\iota_6$), ($\iota_7$) are induced by splittings of
the short exact sequences above and the isomorphism ($\iota_4$) is
the direct sum of the isomorphism $\lambda$ with the identity on
$\oplus_{i \,\, \odd} (B^i_S \oplus B^{i+1}_S) = \oplus_{i \, \,
\even} (B^i_S \oplus B^{i+1}_S)$.

{\bf (3.2) Lemma.} {\em We have:
\[[M^-, \lambda_{M^*}, M^+] = [\Gr(H^-), \lambda,
\Gr(H^+)] \quad {\rm in} \quad G_0(R, \varphi).\] In particular,
the class $[M^-, \lambda_{M^*}, M^+] \in G_0(R,\varphi)$ does not
depend on the chosen splittings.}

{\em Proof.} By relation~(Rii), the class $[M^-, \lambda_{M^*},
M^+]$ can be written  as the sum of the 7 classes $T_i$,
$i=1,\ldots, 7$, of the triples corresponding to the isomorphisms
($\iota_1$), \ldots, ($\iota_7$). We have $T_4 = [\Gr(H^-),
\lambda, \Gr(H^+)]$ by the relation~(Ri) and the relation $[A,
\id, A] =0$. Thus, it suffices to show that $T_1 = T_2 = T_3 = 0 =
T_5 = T_6 = T_7$. This immediately follows from the relation~(Ri)
and the fact that, for any short exact sequence $0 \ra A' \ra A
\ra A'' \ra 0$ of finitely generated $R$-modules, we have $[A,
\alpha, A' \oplus A''] =0$ in $G_0(R,\varphi)$ where the
isomorphism $\alpha: A_S \ra A'_S \oplus A''_S$ is induced by any
splitting of the short exact sequence $0 \ra A'_S \ra A_S \ra
A''_S \ra 0$. To prove this fact, we merely apply relation~(Ri) to
the obvious short exact sequence of triples
\[0 \ra (A',\id, A') \ra (A,\alpha, A'\oplus A'') \ra (A'', \id,
A'') \ra 0\] and use the relations $(A', \id , A') = 0 = (A'',
\id, A'')$. \hfill $\square$

{\bf (3.3) Corollary.} {\em Let $f: \ZZ[G] \hookrightarrow \QQ[G]$
be the canonical inclusion as in Example~(3.1). We assume that
$M^*_\QQ$ is acyclic. Then, under the canonical isomorphism
$G_0(\ZZ[G],\QQ) \cong G_0T(\ZZ[G])$, the class $[M^-, 0_{M^*},
M^+]$ is mapped to $[H^+(M^*)] - [H^-(M^*)]$.}

{\em Proof.} Obvious. \hfill $\square$

We now recall that the {\em relative Grothendieck group
$K_0(R,\varphi)$ of coherent projective $R$-modules} is defined to
be the abelian group with generators $[P,\psi,Q]$ where $P$ and
$Q$ are finitely generated projective $R$-modules and $\psi: P_S
\ra Q_S$ is an $S$-module isomorphism and relations which are
analogous to the relations~(Ri) and (Rii) described above. This
group $K_0(R,\varphi)$ is also often referred to as a relative
algebraic $K$-group.

{\em (3.4) Example.} In the situation of Example~(3.1) we again
write $K_0(\ZZ[G],E)$ for $K_0(\ZZ[G],\varphi)$. It is well-known
and easy to prove that the relative Grothendieck group
$K_0(\ZZ[G],\QQ)$ is isomorphic to the Grothendieck group
$K_0T(\ZZ[G])$ of all finite $\ZZ[G]$-modules of projective
dimension at most 1. Here, the class of a finite $\ZZ[G]$-module
$M$ of projective dimension 1 is mapped to the class of the triple
$(P, \alpha_\QQ, Q)$ where
\[ \xymatrix{0 \ar[r] & P \ar[r]^\alpha & Q \ar[r] & M
\ar[r] & 0}\] is any resolution of $M$ by finitely generated
projective modules. The inverse map sends the class of a triple
$(P,\psi,Q)$ to the element $[\coker(\alpha)] - [P/nP]$ where
$\alpha \in \Hom_{\ZZ[G]}(P,Q)$ and $n \in \ZZ$ are chosen in such
a way that $\psi = \alpha/n$ in $\Hom_{\QQ[G]}(P_\QQ,Q_\QQ) =
\Hom_{\ZZ[G]}(P,Q)_\QQ$. (Again, this statement can be generalized
in the obvious way to the case in which $\ZZ$ is replaced by any
Dedekind domain and $\QQ$ by the corresponding field of
fractions.)

By Theorem~15.5 on p.~216 in \cite{Sw}, we have a natural exact
sequence \begin{equation}\label{lesrkt} K_1(R) \ra K_1(S) \,\,
\stackrel{\partial}{\ra}\,\, K_0(R, \varphi) \ra K_0(R) \ra
K_0(S)\end{equation} where the connecting homomorphism $\partial$
is (uniquely) given by the following rule (see Lemma~15.7 on
p.~217 in \cite{Sw}): Let $P$ be a finitely generated projective
$R$-module and $\psi$ an $S$-module automorphism of $P_S$; then
$\partial$ maps the class of the pair $(P_S,\psi)$ in $K_1(S)$ to
the class of the triple $[P,\psi,P]$ in $K_0(R,\varphi)$.

Let now $P^*$ be a perfect complex of $R$-modules. We assume that,
for each $i$, we are given a (finite, exhaustive and separated)
decreasing filtration $(F^nH^i)_{n \in \ZZ}$ on the cohomology
module $H^i(P^*)$ and that we are also given an $S$-module
isomorphism
\[\lambda: \Gr(H^-)_S \,\, \tilde{\ra} \,\, \Gr(H^+)_S \]
from $\Gr(H^-)_S := \oplus_{i \, \odd} \oplus_{n \in \ZZ}
(F^nH^i/F^{n+1}H^i)_S$ to $\Gr(H^+)_S$. As above, we form the
isomorphism
\[\lambda_{P^*}: P^-_S \,\, \tilde{\ra} \,\, P^+_S\]
which depends on splittings of the corresponding natural short
exact sequences. The following result states that this isomorphism
gives rise to a well-defined class in the relative Grothendieck
Group $K_0(R,\varphi)$ and records the basic properties of this
class. (We remark that this result is an obvious generalization of
  a special case of results proved in \S 1 of \cite{BuIwa}.
However, as that paper is not published, we shall repeat the
relevant arguments here.)

{\bf (3.5) Proposition.} {\em \\ (a) The class
$\chi^\rel(P^*,\lambda) := [P^-, \lambda_{P^*}, P^+] \in
K_0(R,\varphi)$ does not depend upon the chosen splittings.\\ (b)
If the complex $P^*$ is acyclic, then $\chi^\rel(P^*,0) = 0$ in
$K_0(R,\varphi)$.\\ (c) Let $Q^*$ be another perfect complex with
a filtration on the cohomology groups as above and let $\alpha^*:
Q^* \ra P^*$ be a quasi-isomorphism such that $H^i(\alpha^*)$ is
compatible with the given filtrations for all $i$. Then one has
\[[P^-,\lambda_{P^*}, P^-] =
[Q^-, (\Gr(H^+(\alpha^*))^{-1} \circ \lambda \circ
\Gr(H^-(\alpha^*)))_{Q^*}, Q^+]\] in $K_0(R,\varphi)$.
\\ (d) Let $\lambda': \Gr(H^-)_S \,\, \tilde{\ra} \,\,
\Gr(H^+)_S$ be any other $S$-module isomorphism. Then one has
\[\chi^\rel(P^*,\lambda') - \chi^\rel(P^*,\lambda) =
\partial([\Gr(H^-)_S,\lambda^{-1}\circ \lambda']) \quad
\textrm{in} \quad K_0(R,\varphi).\] }

{\em Proof.} \\ (a) Let $0 \ra W' \,\, \stackrel{i}{\ra} \,\, W
\,\, \stackrel{\varepsilon}{\ra} W'' \ra 0$ be a short exact
sequence of finitely generated $S$-modules. Then it is easy to see
that the composition $W \cong W' \oplus W'' \cong W$ of two
isomorphisms induced by any splittings of the given sequence is
given by $w \mapsto w + ih\varepsilon(w)$ for some $S$-module
homomorphism $h: W'' \ra W'$. Therefore, if we choose, for example
for one odd $i$, a different splitting in one of the short exact
sequences
\[ 0 \ra Z^i_S \ra P^i_S \ra B^{i+1}_S \ra 0, \quad
\quad 0 \ra B^i_S \ra Z^i_S \ra H^i_S \ra 0 \] or
\[0 \ra (F^{n+1}H^i)_S \ra (F^nH^i)_S \ra \Gr^n(H^i)_S \ra 0,
\quad n \in \ZZ,\] and denote the corresponding isomorphism by
$\lambda_{P^*}^\dagger$, then we can find a short exact sequence
\[ 0 \ra U \,\, \stackrel{j}{\ra} \,\, P^-_S \,\,
\stackrel{\eta}{\ra} \,\, V \ra 0 \] of $S$-modules such that the
composition $(\lambda_{P^*}^\dagger)^{-1} \lambda_{P^*}$ is given
by $p \mapsto p + j Z \eta(p)$ for some $S$-module homomorphism
$Z: V \ra U$. In particular, we obtain a short exact sequence of
pairs
\[0 \ra (U,\id_U) \ra (P^-_S,(\lambda_{P^*}^\dagger)^{-1}\lambda_{P^*}) \ra
(V,\id_V) \ra 0.\] Hence in $K_0(R,\varphi)$ we have
\begin{eqnarray*}
\lefteqn{[P^-, \lambda_{P^*}, P^+] - [P^-, \lambda_{P^*}^\dagger,
P^+] =[P^-, (\lambda^\dagger_{P^*})^{-1}\circ \lambda_{P^*},P^-]}
\\ &&\hspace*{3em}=\partial([P^-_S, (\lambda_{P^*}^\dagger)^{-1}\circ \lambda_{P^*}])
= \partial([U,\id_U] + [V,\id_V])= 0,
\end{eqnarray*}
as was to be shown.\\
(b) If $P^*$ is acyclic, then for each $i$, the $R$-module $B^i =
Z^i$ is projective and so the short exact sequence of $R$-modules
$0 \ra Z^i \ra P^i \ra B^{i+1} \ra 0$ splits. Using such
splittings we obtain an isomorphism
\[\iota: P^- \cong \oplus_{i \, \odd}(Z^i \oplus B^{i+1}) =
\oplus_{i \, \even} (Z^i \oplus B^{i+1}) \cong P^+\] and, by
claim~(a), we may take $\iota_S$ for $0_{P^*}$. Now, as in the
proof of Lemma~(3.2), we see that $\chi^\rel(P^*,0) = [P^-,
 \iota_S, P^+] = 0$ in $K_0(R,\varphi)$. \\ (c) By Lemma~(3.6)
below, we may assume that, in each degree $i$, both the map
$\alpha^i:Q^i \ra P^i$ and the induced map $Z^i(\alpha^*):
Z^i(Q^*) \ra Z^i(P^*)$ are surjective. Then the kernel $K^*$ of
the epimorphism $\alpha^*:Q^* \ra P^*$ is an acyclic perfect
complex and, in each degree $i$, we have the following commutative
diagram of short exact sequences of $S$-modules
\[\xymatrix{Z^i(K^*)_S \ar@{^{(}->}[r] \ar@{^{(}->}[d]& K^i_S
\ar@{>>}[r] \ar@{^{(}->}[d] & B^{i+1}(K^*)_S \ar@{^{(}->}[d] \\
Z^i(Q^*)_S \ar@{^{(}->}[r] \ar@{>>}[d] & Q^i_S \ar@{>>}[r]
\ar@{>>}[d] & B^{i+1}(Q^*)_S \ar@{>>}[d] \\ Z^i(P^*)_S
\ar@{^{(}->}[r] & P^i_S \ar@{>>}[r] & B^{i+1}(P^*)_S }\] and
\[\xymatrix{B^i(K^*)_S \ar@{^{(}->}[r] \ar@{^{(}->}[d] &
Z^i(K^*)_S \ar@{>>}[r] \ar@{^{(}->}[d] & 0 \ar@{^{(}->}[d] \\
B^i(Q^*)_S \ar@{^{(}->}[r] \ar@{>>}[d] & Z^i(Q^*)_S \ar@{>>}[r]
\ar@{>>}[d] & H^i(Q^*)_S \ar@{>>}[d] \\ B^i(P^*)_S \ar@{^{(}->}[r]
& Z^i(P^*)_S \ar@{>>}[r] & H^i(P^*)_S }.
\] By Lemma~(3.7) below, we can choose compatible splittings in these diagrams. In
particular, we obtain a short exact sequence of triples
\[ \xymatrix{[K^-, 0_{K^*}, K^+] \ar@{^{(}->}[r]&
 [Q^-, (\Gr(H^+(\alpha))^{-1} \circ \lambda \circ
 \Gr(H^-(\alpha)))_{Q^*},Q^+] \ar@{>>}[r] & [P^-,
 \lambda_{P^*},P^+].}
\]
Now, relation~(Ri) and claim~(b) together imply claim~(c).
\\ (d) In $K_0(R,\varphi)$ we have
\[
 \chi^\rel(P^*,\lambda') - \chi^\rel(P^*,\lambda) =  [P^-,
 (\lambda_{P^*})^{-1} \circ \lambda'_{P^*}, P^-]
 = \partial([P^-_S, (\lambda_{P^*})^{-1} \circ \lambda'_{P^*}]).
\]
If we use the same splittings for the definition of
$\lambda'_{P^*}$ and $\lambda_{P^*}$, then the composition of the
isomorphisms ($\iota_1$), ($\iota_2$) and ($\iota_3$) yields an
isomorphism between the pairs $(P^-_S,(\lambda_{P^*})^{-1} \circ
\lambda'_{P^*})$ and
\[(\oplus_{i \, \odd}(\oplus_{n\in \ZZ} (\Gr^n(H^i)_S)
\oplus B^i_S \oplus B^{i+1}_S),\lambda^{-1} \circ \lambda'
\oplus\id_{\oplus_{i \, \odd} (B^i_S \oplus B^{i+1}_S)}).\] Hence,
we have $[P^-_S, (\lambda_{P^*})^{-1} \circ \lambda'_{P^*}] =
[\Gr(H^-)_S,\lambda^{-1}\circ\lambda']$ in $K_1(S)$. This
completes the proof of claim~(d) of Proposition~(3.5). \hfill
$\square$

{\bf (3.6) Lemma.} {\em Let $R$ be a ring and $\alpha^* :Q^* \ra
P^*$ a quasi-isomorphism between perfect complexes of $R$-modules.
Then there exists a perfect complex $T^*$ of $R$-modules together
with quasi-isomorphisms $\beta^* : T^* \ra Q^*$ and $\gamma^*: T^*
\ra P^*$ such that, in each degree $i$, we have $H^i(\alpha^*
\circ \beta^*) = H^i(\gamma^*)$ and both the maps $\beta^i:T^i \ra
Q^i$, $Z^i(\beta^*): Z^i(T^*) \ra Z^i(Q^*)$ and the maps
$\gamma^i: T^i \ra P^i$, $Z^i(\gamma^*): Z^i(T^*) \ra Z^i(P^*)$
are surjective.}

{\em Proof.} Let $K^*$ be the acyclic perfect complex with
$K^i:=P^{i-1} \oplus P^i$ and with differential $K^i \ra K^{i+1}$,
$(x_{i-1},x_i) \mapsto (x_i,0)$. We set $T^*:= K^* \oplus Q^*$ and
define $\beta^*: T^* \ra Q^*$ to be the canonical projection. Then
$\beta^*$ is clearly a quasi-isomorphism and, for each $i$, the
maps $\beta^i$ and $Z^i(\beta^*)$ are surjective. We define
$\gamma^i: T^i=P^{i-1} \oplus P^i \oplus Q^i \ra P^i$ by
$\gamma^i|_{P^{i-1}} = d^{i-1}(P^*)$, $\gamma^i|_{P^i} =
\id_{P^i}$, and $\gamma^i|_{Q^i} = \alpha^i$. Then, $\gamma^*$ is
clearly a homomorphism of complexes such that, in each degree $i$,
we have $H^i(\alpha^* \circ \beta^*) = H^i(\gamma^*)$ and the map
$\gamma^i$ is surjective. One easily checks that the map
$Z^i(\gamma^*): P^{i-1} \oplus Z^i(Q^*) \ra Z^i(P^*)$ is also
surjective. Thus, the proof of Lemma~(3.6) is complete. \hfill
$\square$

{\bf (3.7) Lemma.} {\em Let $S$ be a semisimple ring. Suppose we
are given a commutative diagram of short exact sequences of
$S$-modules
\[\xymatrix{K' \ar@{^{(}->}[r] \ar@{^{(}->}[d]& K
\ar@{>>}[r]^\delta \ar@{^{(}->}[d] & K'' \ar@{^{(}->}[d] \\ V'
\ar@{^{(}->}[r] \ar@{>>}[d] & V \ar@{>>}[r]^\varepsilon
\ar@{>>}[d] & V'' \ar@{>>}[d] \\ W' \ar@{^{(}->}[r] & W
\ar@{>>}[r]^\eta & W''. }\] Then one can choose a section $\sigma:
V'' \ra V$ to $\varepsilon$ such that $\sigma|_{K''}$ is a section
to $\delta$ (and hence $\sigma$ induces a section to $\eta$).}

{\em Proof.} First, we choose any section $\tilde{\sigma}$ to
$\varepsilon$. Then the composition
\[\xymatrix{K''
\ar@{^{(}->}[r] & V'' \ar[r]^{\tilde{\sigma}} & V \ar@{>>}[r] &
W}\] may be considered as a homomorphism $\theta$ from $K''$ to
$W'$. Since $K''$ is projective, we may lift $\theta$ to a
homomorphism from $K''$ to $V'$ and, since $V'$ is injective, this
can in turn be extended to a homomorphism from $V''$ to $V'$ which
we view as a homomorphism $\tilde{\theta}$ from $V''$ to $V$. Now,
one easily checks that the homomorphism $\sigma := \tilde{\sigma}
- \tilde{\theta}$ is a section to $\varepsilon$ such that
$\sigma|_{K''}$ is a section to $\delta$. \hfill $\square$

We now return to the context described at the beginning of this
section. We thus assume given an $E$-trivialized nearly perfect
complex $(C^*,\lambda)$ with $C^* = (C^*, (L_i)_i, (\tau_i)_i)$.
 We choose a quasi-isomorphism $P^* \ra C^*$ as in
 \S2. From Proposition~(2.2) we then obtain natural exact
sequences
\begin{equation}\label{s1} 0 \ra H^i(C^*)_\codiv \otimes \ZZ_l
\ra H^i(\hat{P}^*) \ra \Hom(L_{i+1},\ZZ_l) \ra 0, \quad i \in \ZZ.
\end{equation}
Furthermore, we choose resolutions $\xymatrix{0 \ar[r]& R^{i-1}
\ar[r] & Q^i \ar@{>>}[r]^>>>>{\varepsilon^i} & \Hom(L_i,\QQ)
\ar[r] & 0}$, $i\in \ZZ$, and a homomorphism of complexes
$(\alpha^i,\beta^i)_i: (Q\oplus R)^* \ra P^*$ such that the
composition with the chosen quasi-isomorphism $P^* \ra C^*$ is of
the form considered in \S1. By tensoring the exact sequences~(1)
in \S1 with $\ZZ_l$, we thus obtain natural exact sequences
\begin{equation}\label{s2}
0 \ra H^i(C^*)_\codiv \otimes \ZZ_l \ra H^i(\Cone_C^*) \otimes
\ZZ_l \ra \Hom(L_{i+1}, \ZZ_l) \ra 0, \quad i \in \ZZ.
\end{equation}

{\bf (3.8) Proposition.} {\em For each $i$, the
extension~(\ref{s1}) is the negative of the extension~(\ref{s2})
in
\[\Ext^1_{\ZZ_l[G]}(\Hom(L_{i+1},\ZZ_l),
H^i(C^*)_\codiv \otimes \ZZ_l).\]}

{\em Proof.} Let $\Cone_P^*$ denote the mapping cone of the map
$(\alpha^i,\beta^i)_i: (Q\oplus R)^* \ra P^*$. It obviously
suffices to show that the extension
\begin{equation}\label{s3}
0 \ra H^i(P^*)_\codiv \otimes \ZZ_l \ra H^i(\hat{P}^*) \ra
\Hom(L_{i+1},\ZZ_l) \ra 0
\end{equation}
(resulting from Proposition~(2.2)) is the negative of the
extension
\begin{equation}\label{s4}
0 \ra H^i(P^*)_\codiv \otimes \ZZ_l \ra H^i(\Cone_P^*) \otimes
\ZZ_l \ra \Hom(L_{i+1},\ZZ_l) \ra 0
\end{equation}
(constructed as in \S1). We have the following commutative diagram
\[\xymatrix{\widehat{H^i(P^*)_\codiv} \ar@{^{(}->}[r]
\ar@{^{(}->}[d] & H^i(\hat{P}^*) \ar@{>>}[r] \ar[d]^\wr &
\Hom(L_{i+1},\ZZ_l) \\ \widehat{H^i(\Cone_P^*)} \ar[r]^\sim
\ar@{>>}[d] & H^i(\widehat{\Cone}_P^*) \\ \Hom(L_{i+1},\ZZ_l) }\]
where the top row is the extension (\ref{s3}), the left column is
the extension (\ref{s4}), the isomorphism in the second row is a
consequence of Proposition~(2.2) and the isomorphism in the second
column is induced by the canonical inclusion $P^* \hookrightarrow
\Cone_P^*$ (see the proof of Theorem~(2.4)). The snake lemma
applied to this diagram yields a $\ZZ _l[G]$-module automorphism
of $\Hom(L_{i+1},\ZZ_l)$. It suffices to show that this
automorphism is the multiplication by $-1$. For this it suffices
to check that the connecting homomorphism from the upper right
corner in the diagram
\[\xymatrix{&&\Hom(L_{i+1}, \ZZ/l^n\ZZ) \ar@{^{(}->}[d]^{\tau_{i+1}} \\
H^i(P^*)/l^nH^i(P^*) \ar@{^{(}->}[r] \ar@{^{(}->}[d] &
H^i(P^*/l^nP^*) \ar@{>>}[r]^\partial \ar[d]^\wr &
{}_{l^n}H^{i+1}(P^*) \ar@{>>}[d] \\
H^i(\Cone_P^*)/l^nH^i(\Cone_P^*) \ar@{^{(}->}[r] \ar@{>>}[d] &
H^i(\Cone_P^*/l^n\Cone_P^*) \ar@{>>}[r]^\partial &
{}_{l^n}H^{i+1}(\Cone_P^*)  \\ \Hom(L_{i+1},\ZZ/l^n\ZZ) }\] (with
the obvious maps) to the lower left corner is the multiplication with $-1$.\\
Any element of $\Hom(L_{i+1},\ZZ/l^n\ZZ)$ can be written as the
residue class $\bar{e}$ of some $e \in \Hom(L_{i+1},\ZZ) \subseteq
\Hom(L_{i+1},\QQ)$. We choose an element $q \in Q^{i+1}$ with
$\varepsilon^{i+1}(q) = e /l^n$ in $\Hom(L_{i+1},\QQ)$. Under the
inclusion $\tau_{i+1}: \Hom(L_{i+1},\ZZ/l^n\ZZ) \hookrightarrow
{}_{l^n}H^{i+1}(P^*)$, the element $\bar{e}$ is then mapped to
$\tau_{i+1}(\overline{\varepsilon^{i+1}(q)})$ where
$\overline{\varepsilon^{i+1}(q)}$ denotes the image of
$\varepsilon^{i+1}(q)$ in $\Hom(L_{i+1},\QQ/\ZZ)$. Since $l^n
\cdot \tau_{i+1}(\overline{\varepsilon^{i+1}(q)}) = 0$ in
$H^{i+1}(P^*)$, we can find an element $p \in P^i$ such that $l^n
\cdot \alpha^{i+1}(q) = d(p)$ in $P^{i+1}$. Then the cohomology
class $[\bar{p}] \in H^i(P^*/l^nP^*)$ of the cocycle $\bar{p} \in
P^i/l^nP^i$ is mapped to
$\tau_{i+1}(\overline{\varepsilon^{i+1}(q)})$ under the connecting
homomorphism $\partial$ in the upper row. Furthermore, the
cohomology class of the cocycle $(p,-l^n q, 0) \in P^i \oplus
Q^{i+1} \oplus R^{i+1} = \Cone_P^i$ is obviously mapped to the
cohomology class of $(\bar{p}, 0, 0)$ under the canonical
inclusion
\[H^i(\Cone_P^*)/l^nH^i(\Cone_P^*) \hookrightarrow
H^i(\Cone_P^*/l^n\Cone_P^*).\] Finally, the cohomology class of
$(p, -l^nq,0)$ is mapped to $-l^n\varepsilon^{i+1}(q) = -e$ under
the epimorphism $\xymatrix{H^i(\Cone_P^*) \ar@{>>}[r] &
\Hom(L_{i+1},\ZZ)}$. Therefore, $\bar{e} \in
\Hom(L_{i+1},\ZZ/l^n\ZZ)$ is mapped to $-\bar{e}$ under the
connecting homomorphism in the diagram above, as was to be shown.
\hfill $\square$

We next observe that, by the classical Noether-Deuring theorem
(see Theorem~(29.7) on p.~200 in \cite{CR}), the
 existence of an isomorphism $\lambda$ as in (\ref{Etriv}) implies
 that there also exists a (non-canonical) $\QQ[G]$-equivariant
isomorphism
\begin{equation}\label{Qtriv} \tilde{\lambda}:\QQ \otimes H^-(C^*)_\codiv \oplus
\Hom(L_+,\QQ) \,\, \tilde{\ra} \,\, \QQ \otimes H^+(C^*)_\codiv
\oplus \Hom(L_-,\QQ).\end{equation} For each prime $l$ we may
therefore tensor $\tilde{\lambda}$ with $\QQ_l$ in order to obtain
a $\QQ_l[G]$-equivariant isomorphism
\[\tilde\lambda_l: \QQ_l \otimes H^-(C^*)_\codiv
\oplus \Hom(L_+,\QQ_l) \,\, \tilde{\ra} \,\, \QQ_l \otimes
H^+(C^*)_\codiv \oplus \Hom(L_-,\QQ_l).\]

As in \S2, we choose a quasi-isomorphism $\tilde{P}^* \ra
\hat{P}^*$ from a perfect complex $\tilde{P}^*$ of
$\ZZ_l[G]$-modules to the $l$-adically completed complex
$\hat{P}^*$. For each $i$, we consider the extension (\ref{s1}) as
a 2-step filtration on $H^i(\tilde{P}^*)$. Then, from the
construction above, we obtain a well-defined element
\[\chi^\rel(C^*,\tilde{\lambda}_l) := [\tilde{P}^-, (\tilde\lambda_l)_{\tilde P^*},
\tilde{P}^+] \in K_0(\ZZ_l[G], \QQ_l)\] which, by Lemma~(1.2),
Corollary~(2.3) and Proposition~(3.5)(c), does not depend on the
chosen quasi-isomorphisms $P^* \ra C^*$ and $\tilde{P}^* \ra
\hat{P}^*$. Finally we recall that there are canonical
isomorphisms
\[\oplusm_{l\, {\rm prime}} K_0(\ZZ_l[G],\QQ_l) \cong \oplusm_{l \,
{\rm prime}} K_0T(\ZZ_l[G]) \cong K_0T(\ZZ[G]) \cong
K_0(\ZZ[G],\QQ),\] and a canonical injective homomorphism
$i_E:K_0(\ZZ[G],\QQ) \ra K_0(\ZZ[G],E)$.

As promised at the beginning of this section, we now define a
canonical refined Euler characteristic $\chi^{\rel}(C^*,\lambda)$
 for the $E$-trivialized nearly perfect complex $(C^*,\lambda)$.

{\bf (3.9) Theorem.} \\ (a) {\em The tuple
\[\chi^\rel(C^*,\tilde\lambda) := (\chi^\rel(C^*,\tilde\lambda_l))_{l \,
{\rm prime}} \in \prod_{l \, {\rm prime}} K_0(\ZZ_l[G],\QQ_l)\]
belongs to the direct sum $\oplus_{l \, {\rm prime}}
K_0(\ZZ_l[G],\QQ_l) \cong K_0(\ZZ[G],\QQ)$.}
\\ (b) {\em The element
\[\chi^\rel(C^*,\lambda):= i_E(\chi^\rel(C^*,\tilde{\lambda})) +
\partial([E \otimes H^-(C^*)_\codiv \oplus \Hom(L_+,E),
\tilde{\lambda}_E^{-1} \circ {\lambda}])\] of $K_0(\ZZ[G],E)$
depends only upon $C^*$ and $\lambda$.}
\\ (c) {\em The image of $\chi^\rel(C^*,\lambda)$ under the
canonical map $K_0(\ZZ[G],\QQ) \ra K_0(\ZZ[G])$, $[P,\psi,Q]
\mapsto [Q] - [P]$, is equal to the Euler characteristic
$\chi(C^*)$.}
\\ (d) {\em The image of $\chi^\rel(C^*,\lambda)$ under the
forgetful map $K_0(\ZZ[G],E) \ra G_0(\ZZ[G],E)$ is equal to the
element
\[[H^-(C^*)_\codiv \oplus \Hom(L_+,\ZZ), \lambda,
H^+(C^*)_\codiv \oplus \Hom(L_-,\ZZ)].\]}

{\em Proof.} As in the proof of Theorem~(2.4), we choose a
quasi-isomorphism $P_\Cone^* \ra \Cone_P^*$ from a perfect complex
$P_\Cone^*$ of $\ZZ[G]$-modules to the mapping cone $\Cone_P^*$ of
the homomorphism $(\alpha^i,\beta^i)_i: (Q\oplus R)^* \ra P^*$
chosen above. We recall from \S1 that we have short exact
sequences
\[0 \ra H^i(C^*)_\codiv \ra H^i(P_\Cone^*) \ra \Hom(L_{i+1},\ZZ)
\ra 0, \quad i \in \ZZ.\] From the construction above, we
therefore obtain an element
\begin{equation}\label{element}[P_\Cone^-, \tilde\lambda_{P^*_{\rm
Cone}}, P_\Cone^+] \in K_0(\ZZ[G],\QQ).\end{equation} To prove
claim~(a) it therefore suffices to show that this element is equal
to $\chi^\rel(C^*,\tilde\lambda)$. In other words, it suffices to
show that for all primes $l$ the element
\[[P_\Cone^- \otimes \ZZ_l, (\tilde\lambda_{P_{\rm Cone}^*})_{\QQ_l},
P_\Cone^+ \otimes \ZZ_l] = [\widehat{P_\Cone}^-,
(\tilde{\lambda}_l)_{\widehat {P_{\rm Cone}}^*},
\widehat{P_\Cone}^+]\] coincides with
$\chi^\rel(C^*,\tilde\lambda_l)$. But the natural inclusion $P^*
\hookrightarrow \Cone_P^*$ induces a quasi-isomorphism $\hat{P}^*
\ra \widehat{\Cone_P}^*$ (see the proof of Theorem~(2.4)). So, by
Lemma~(1.2), we obtain a quasi-isomorphism between the perfect
complexes $\tilde{P}^*$ and $\widehat{P_\Cone}^*$ which, by (the
proof of) Proposition~(3.8), is compatible with the filtrations on
the cohomology given by the short exact sequences (\ref{s1}) and
(\ref{s2}). Thus the desired equality follows from
Proposition~(3.5)(c).\\
We now turn to claim~(b). Here it suffices to show that the
element $\chi^\rel(C^*,\lambda)$ is independent of the choice of
$\QQ [G]$-equivariant isomorphism $\tilde \lambda$ as in
(\ref{Qtriv}). But if $\lambda^\dagger$ is any other isomorphism
as in (\ref{Qtriv}), then Proposition~(3.5)(d) implies that
\[i_E(\chi^\rel(C^*,\lambda^\dagger)) = i_E(\chi^\rel(C^*,\tilde{\lambda})) +
\partial([E \otimes H^-(C^*)_\codiv \oplus \Hom(L_+,E),
\tilde{\lambda}_E^{-1} \circ \lambda_E^\dagger])\] in
$K_0(\ZZ[G],E)$. This leads directly to the required equality
 since in $K_1(E[G])$ one has
\[[E \otimes H^-(C^*)_\codiv \oplus \Hom(L_+,E),
(\lambda^\dagger_E)^{-1} \circ \lambda] + [E \otimes
H^-(C^*)_\codiv \oplus \Hom(L_+,E), \tilde{\lambda}_E^{-1} \circ
\lambda_E^\dagger]\] \[ = [E \otimes H^-(C^*)_\codiv \oplus
\Hom(L_+,E), \tilde{\lambda}_E^{-1} \circ \lambda].\] Next we
observe that the exact sequence of relative $K$-theory
(\ref{lesrkt}) implies that the elements $\chi^\rel(C^*,\lambda)$
and $\chi^\rel(C^*,\tilde{\lambda})$ have the same image under the
 homomorphism $K_0(\ZZ [G],E) \ra K_0(\ZZ [G])$. Claim~(c)
therefore follows from the fact that $\chi^{\rm rel}
(C^*,\tilde{\lambda})$ is equal to the element which occurs in
(\ref{element}). Indeed, the latter element clearly maps to
 $\chi(P_\Cone^*)$ in $K_0(\ZZ[G])$ which is equal to $\chi(C^*)$
by Theorem~(1.3).\\
To prove claim~(d), we first observe that the element
({\ref{element}) maps to the element
$[P^-_\Cone,\tilde\lambda_{P_{\rm Cone}},P^+_\Cone]$ in
$G_0(\ZZ[G],\QQ)$. Further, Lemma~(3.2) implies that the latter
element is equal to
\begin{eqnarray*}
\lefteqn{[\Gr(H^-(\Cone_C^*)),\tilde\lambda,\Gr(H^+(\Cone_C^*))] }\\
&=&[H^-(C^*)_\codiv\oplus
\Hom(L_+,\ZZ),\tilde\lambda,H^+(C^*)_\codiv \oplus
\Hom(L_-,\ZZ)].\end{eqnarray*}

It therefore follows that the image of $\chi^{\rm
rel}(C^*,\lambda)$ in $G_0(\ZZ [G],E)$ is equal to
\begin{eqnarray*}
\lefteqn{[H^-(C^*)_\codiv \oplus
\Hom(L_+,\ZZ),\tilde\lambda_E,H^+(C^*)_\codiv \oplus
\Hom(L_-,\ZZ)]}\\ &&+ [H^-(C^*)_\codiv \oplus \Hom(L_+,\ZZ),
\tilde{\lambda}_E^{-1}\circ{\lambda}, H^-(C^*)_\codiv \oplus
\Hom(L_+,\ZZ)]\\ &= &[H^-(C^*)_\codiv \oplus \Hom(L_+,\ZZ),
\lambda, H^+(C^*)_\codiv \oplus \Hom(L_-,\ZZ)],\end{eqnarray*} as
required.\\
This completes the proof of Theorem~(3.9). \hfill $\square$

\bigskip

\bigskip

David Burns, Department of Mathematics, King's College London,
Strand, London WC2R 2LS, United Kingdom. \\
{\em E-mail:} david.burns@kcl.ac.uk

Bernhard K\"ock, Faculty of Mathematical Studies, University of
Southampton, Highfield, Southampton SO17 1BJ, United Kingdom.\\
{\em E-mail:} bk@maths.soton.ac.uk

Victor Snaith, Faculty of Mathematical Studies, University of
Southampton, Highfield, Southampton SO17 1BJ, United Kingdom.\\
{\em E-mail:} vps@maths.soton.ac.uk

\end{document}